\documentclass[aps,pra,superscriptaddress,nofootinbib]{revtex4}
\pdfoutput=1
\usepackage[intlimits]{amsmath}
\usepackage{amsfonts}

\usepackage{enumitem}
\usepackage{psfrag}

\usepackage{amsthm}
\usepackage{float}
\usepackage{tikz}

\advance\voffset by -1cm
\advance\hoffset by -0.5cm
\newcommand\beq            {\begin{equation}}
\newcommand\be            {\begin{equation}}
\newcommand\bea           {\begin{equation}\begin{array}l\displaystyle}
\newcommand\ee            {\end{equation}}
\newcommand\eeq            {\end{equation}}
\newcommand\bes           {\begin{subequations}}
\newcommand\esu           {\end{subequations}}

\renewcommand{\(}{\left(}
\renewcommand{\)}{\right)}
\renewcommand{\[}{\left[}
\renewcommand{\]}{\right]}

\newcommand{\bigx}[1]{\bBigg@{#1}}

\def\3pt#1#2#3{{\langle{#1}\vert{#2}\vert{#3}\rangle}}

\newcommand\doi[2]        {\href{http://dx.doi.org/#1}{#2}}

\newcommand{\EQ}{\begin{equation}}
\newcommand{\EN}{\end{equation}}
\usepackage{epsfig}
\usepackage{psfrag}
\usepackage{amsmath}
\usepackage{graphicx}
\usepackage{amsfonts}
\usepackage{amssymb}

\def\tilde{\widetilde}

\def\hat{\widehat}
\def\*{\star}
\def\[{\left[}
\def\]{\right]}
\def\({\left(}      

\def\){\right)}

\def\frac#1#2{\dfrac{#1}{#2}}
\def\inv#1{\dfrac{1}{#1}}
\def\half{\tfrac{1}{2}}

\def\2pi{\hbox{$2\pi i$}}

\def\dsl{\raise.15ex\hbox{/}\kern-.57em\partial}
\def\Dsl{\,\raise.15ex\hbox{/}\mkern-.13.5mu D}
\def\be{\beta}

   \def\CE{{\cal E}}

   \def\CN{{\cal N}}

\def\2pi{\hbox{$2\pi i$}}

\def\dsl{\raise.15ex\hbox{/}\kern-.57em\partial}
\def\Dsl{\,\raise.15ex\hbox{/}\mkern-.13.5mu D}

\def\barray{\begin{eqnarray}}
\def\earray{\end{eqnarray}}
\def\beq{\begin{equation}}
\def\eeq{\end{equation}}

\def\AA{\leavevmode\setbox0=\hbox{h}
\dimen0=\ht0 \advance\dimen0 by-1ex\rlap{\raise.67\dimen0\hbox{\char'27}}A}

\def\half{\tfrac{1}{2}}

\def\prob{{ \bf Pr}}
\def\Ex{{\bf E}}

\def\half{\tfrac{1}{2}}

\def\dist{{\overset{d}{\longrightarrow}}}

\def\s{\ell} 

\newtheorem{theorem}{Theorem}

\newtheorem{conjecture}{Conjecture}

\begin{document}
\bibliographystyle{plainnat}

\title{{\Large {\bf Generalized Riemann Hypothesis  
 \\
and Stochastic Time Series}}} 

\author{Giuseppe Mussardo}
\affiliation{SISSA and INFN, Sezione di Trieste, via Bonomea 265, I-34136, 
Trieste, Italy}
\affiliation{Institute for Theoretical Physics, Utrecht University}
\author{ Andr\'e  LeClair}
\affiliation{Cornell University, Physics Department, Ithaca, NY 14850}

\begin{abstract}
\noindent
\noindent
Using the Dirichlet theorem on the equidistribution of residue classes modulo $q$ and the Lemke Oliver-Soundararajan conjecture on the distribution of pairs of residues on consecutive primes, we show that the domain of convergence of the infinite product of Dirichlet $L$-functions of non-principal characters can be extended from $\Re(s) > 1$ down to $\Re(s)  > \half$, without encountering any zeros before reaching this critical line. The possibility of doing so can be traced back to a universal diffusive random walk behavior  $C_N  = {\cal O}(N^{1/2})$ of  the series $C_N = \sum_{n=1}^N \chi (p_n)$ over the primes $p_n$ where $\chi$ is a Dirichlet character,  which underlies the convergence of the infinite product of the Dirichlet functions. The series $C_N$ presents several aspects in common with stochastic time series and its control requires to address a problem similar to the Single Brownian Trajectory Problem in statistical mechanics. In the case of the Dirichlet functions of non principal characters, we show that this problem can be solved in terms of a self-averaging procedure based on an ensemble $\CE$ of block variables computed on extended intervals of primes. Those intervals, called {\em inertial intervals}, ensure the ergodicity and stationarity of the time series underlying the quantity $C_N$.  The infinity of primes also ensures the absence of rare events which would have been responsible for a different scaling behavior than the universal law $C_N  = {\cal O}(N^{1/2})$ of the random walks.

\end{abstract}

\maketitle

\section{Introduction}

The Generalised Riemann Hypothesis (GRH) states that the non-trivial zeros of {\em all} the Dirichlet 
$L$-functions $L(s,\chi)$ of the complex variable $s = \sigma + i t$ lie along the critical line $\Re(s) 
= \half$. In this paper we focus our attention on $L$-functions associated to non-principal characters (see 
Appendix \ref{Appcharacters} for their definitions). As we are going to show below, the fact that the non-trivial zeros of these functions align themselves along such a critical line finds a natural and universal explanation in terms of a random walk, in the sense that $\sigma =\half$ is the universal critical exponent of a random walk process $C_N = {\cal O}(N^{1/2})$ which exists for any Dirichlet $L$-function of a  non principal character. Such a scaling behavior of the quantity $C_N$ is a consequence of the Dirichlet theorem on the equi-distribution of reduced residue classes modulo $q$ \cite{Diric,SelbergD} and the Lemke Oliver-Soundararajan  (LOS) conjecture on the distribution of pairs of residues on consecutive primes \cite{OliverSoundararajan}. Moreover, there are very intriguing analogies between the deterministic series $C_N$ coming from the Dirichlet $L$-functions and the stochastic time series. The close relation between these two entities will play a crucial role in the analysis that follows. From a statistical mechanics point of view, there is also a very close and interesting parallel of the series $C_N$ with the Single Brownian Trajectory Problem (see, for istance \cite{brow1,brow2,brow3,brow4} and references therein), problem which consists in extracting the statistical properties of a sequence of events which happen only once!

This paper is devoted to a short presentation of our findings,  whereas  in  \cite{longpaper} we provide a thorough discussion of many other aspects of $L$-functions,  a detailed derivation of the results discussed here, along with extensive numerical evidence.   For general definitions and properties of the Dirichlet $L$ functions, we refer the reader to standard references, for instance \cite{Apostol,Iwaniec,Bombieri,Steuding,Sarnak}. 

\vspace{3mm}
\noindent
{\bf $L$-functions}. For  $\Re(s) > 1$ there are two equivalent representations of these functions, one given in terms of an infinite series on the natural numbers $m$, the other in terms of an infinite product over the sequence of primes $p_n$ (hereafter labelled in ascending order) 
\beq
L(s,\chi) \,=\,
\sum_{m=1}^\infty \frac{\chi(m)}{m^s}\,=\, 
 \prod_{n=1}^\infty  \( 1 -  \frac{\chi (p_n)}{p_n^s} \)^{-1} \,\,\,. 
\label{Euleridentity}
\eeq
This equation expresses the so-called {\em Euler identity}. The quantities $\chi(m)$ which enter the definition of $L(s,\chi)$ are  the so-called Dirichlet characters. When non-zero, they are pure phases and therefore expressed in terms of some angles $\theta_m$ defined as 
\beq
\label{thetan}
\chi (m) \,=\,  e^{i \theta_m }, ~~~~~\forall ~\chi(m) \neq 0 \,\,\,.
\eeq   
The characters $\chi(m)$ depend on a positive integer $q$, called the {\em modulus}. For any given $q$, there are $\varphi(q)$ different characters which correspond to the irreducible representation of the abelian group given by the prime residue classes modulo $q$    
\beq
(\mathbb{Z}/q\mathbb{Z})^* := \{m\, {\rm mod} \,q \,:\, (m, q) = 1\} \,\,\,, 
\label{groupab}
\eeq
where $\varphi(q)$ -- the Euler totient arithmetic function --  is the dimension of this group. One can focus only on the cases where $q$ is a prime number, where  $\varphi(q) = q-1$, since it can be shown that the main properties of all other $L$-functions can be derived from  these cases (when $q$ is prime, the characters are primitive characters). Dirichlet $L$-functions of non-principal characters are entire functions (they have only zeros in the complex plane $s$) while those associated to the  principal character at any $q$, in addition to zeros,  also have  a pole at $s=1$. Notice that the Riemann $\zeta$ function \cite{Riemann1,Riemann2,Riemann3} is a particular Dirichlet $L$-function, the one associated to the principal character with  $q=1$. Other important properties both of the characters and Dirichlet $L$-functions are collected in the various appendices at the end of this  paper. 

\vspace{3mm}
\noindent
{\bf Zeros and spectra}. 
Many authors have studied in detail the zeros of  Dirichlet $L$-functions, obtaining precise  information on their locations: Selberg \cite{Selberg1} and then Fuji \cite{Fujii} obtained the counting formula $N(T,\chi)$ for the number of zeros up to height $T$ within the entire critical strip $0 \le \Re(s) \le 1$; Iwaniec, Luo and Sarnak \cite{IwaniecS} studied the distribution of low lying zeros of $L$-functions near and at the critical line; Conrey and, in a separate paper, Hughes and Rudnick \cite{Conrey,Hughes}, following the original Montgomery-Odlyzko conjecture relative to the zeros of the Riemann $\zeta$-function  and their relation to random matrix theory \cite{Dyson,Montgomery,Odlyzko} (see also \cite{Rudnick-Sarnack}), studied the statistics of the zeros of the $L$-function; Conrey, Iwaniec, and Soundararajan estimated that more than $56\%$  of the non-trivial zeros are on the critical line \cite{Conrey2}. An exact transcendental equation for the $n$-th zero,  which depends only on $n$,  was obtained in \cite{Transcendental}. This, of course, is not at all an exhaustive list of all  the authors who have contributed to this subject and we apologize for those not cited here.   

Apart from analysis,  another line of attack to the problem, although mostly focused on the Riemann $\zeta$-function, has been based on the original suggestion by P\'olya \cite{Polya},  and concerns the deep interplay between the spectral theory of quantum mechanics and the zeros of the Riemann $\zeta$ function. This viewpoint has given rise to an important series of works by Berry, Keating, Connes, Sierra, Srednicki, Bender and many others \cite{BK1,BK2,BK3,BK4,Connes,Sierra1,Sierra2,Srednicki,Bender} (for a more complete list of references, see the review \cite{reviewRiemann}). In these approaches, the main focus consists of identifying a quantum mechanical Hamiltonian whose spectrum coincides with the imaginary parts of the Riemann zeros along the line $\Re(s) = \half$, in order to argue that such a distribution of the Riemann zeros is a consequence of the spectral theory of quantum Hamiltonians. The approach presented in this paper, however, is of a completely different nature,  since it rather has strong analogies in statistical mechanics.

\vspace{3mm}
\noindent
{\bf Main idea of the paper}. In this paper we are concerned with Dirichlet $L$-functions relative to non-principal characters (the technical reason behind this choice will become clear below). Our line of attack on the GRH for these functions consists of the possibility to extend the convergence of their infinite product representation for $\Re(s) < 1 $ and show that this can be done down to the critical line $\Re(s) =\half$,  without encountering  any zeros along the way. In view of the functional equation (\ref{FELambda}) which links $L(s,\chi)$ to $L(1-s,\overline\chi)$, this implies that all non-trivial zeros of the Dirichlet $L$-functions relative to non-principal characters are on the critical line $\Re(s) =\half$. 

\vspace{3mm}
\noindent
{\bf Organization of the paper}. The paper is organized as follows. In Section \ref{extending} we discuss the 
possibility to extend the infinite product representation of the Dirichlet $L$-function to the left of the half-plane 
$\Re(s) =1$, identifying the quantity $C_N$, defined in eq.\,(\ref{Gseries}), which is a direct diagnosis for the real part of the first zero encountered in the critical strip. The quantity $C_N$ depends on the phases $\theta_{p_n}$ of the characters computed on primes and Section \ref{stprth} is devoted to the statistical analysis of these quantities: we recall, in particular, the Dirichlet theorem on the equi-distribution of reduced residue classes modulo $q$ and the Lemke Oliver-Soundararajan conjecture on the distribution of pairs of residues on consecutive primes. Section \ref{Timeseries} is the key part of this paper: after a brief discussion of the Single Brownian Trajectory Problem, we define the block variables $C_N(l)$ and from  them  define a statistical ensemble $\CE$ for the original quantity $C_N$. The key properties of this ensemble is the zero mean of 
$C_N$ and a {\em linear} behaviour in $N$ (up to logarithmic corrections) of the variance of this quantity. 
This implies that, for large $N$, the quantity $C_N$ goes as $C_N = {\mathcal O}(N^{1/2})$  and therefore that the domain of convergence of the infinite product of Dirichlet $L$-functions of non-principal characters can be extended from $\Re(s) > 1$ down to $\Re(s)  > \half$, without encountering any zeros before reaching this critical line. The paper has also three appendices devoted  to various aspects of the Dirichlet characters and Dirichlet $L$-functions.

\section{Extending the infinite product into the critical strip}\label{extending}
\label{Probabilistic}
In this Section, following the reference \cite{EPFchi}, we identify a quantity -- in particular a series --   whose proposed properties allow for the extension of convergence of the Euler product of the $L$-functions into the critical strip,  and is therefore an ideal marker for the location of their zeros. 

Consider the infinite product representation of the $L$-functions 
\beq 
L(s,\chi) \,=\, 
 \prod_{n=1}^\infty  \( 1 -  \frac{\chi (p_n)}{p_n^s} \)^{-1} \,\,\,,
\eeq
and take the logarithm on both sides of this equation, so that 
\beq 
\log L(s,\chi) \,=\,  X(s,\chi) + R(s,\chi)\,\,\,,
\eeq
where
\beq\label{PDir}
X(s,\chi) = \sum_{n=1}^{\infty} \dfrac{\chi(p_n)}{p_n^{\, s}}\,\,\,\,\,\,\, , \qquad
R(s,\chi) = \sum_{n=1}^{\infty} \sum_{m=2}^{\infty}
\dfrac{\chi(p_n)^{m}}{m p_n^{\, ms}}\,\,\,.
\eeq
Since $R(s,\chi)$ absolutely converges for $\sigma > \half$, we have 
\beq\label{logsum2}
\log L(s,\chi) = X(s,\chi) + O(1),
\eeq
namely, the convergence of the Euler product to the right of the critical line depends only on $X(s,\chi)$. On the other hand, the singularities of $ \log L(s,\chi)$ are determined by the zeros of $L(s,\chi)$,  and, if present, also by the pole $s=1$. Concerning  the GRH, the main emphasis is of course in locating the {\em zeros} of these functions and therefore the presence of a pole at $s=1$ just introduces a simple, though significant, complication\footnote{See the discussion below, after eq.\,(\ref{noshiftfunction}), for the effect induced by the pole at $s=1$. This affects the discussion of the Riemann Hypothesis for the original Riemann $\zeta$ function.}. The advantage of the $L$-functions of non-principal characters is that they do not have a pole at $s=1$ (see Appendix C) and therefore for all of them we have a very succinct mathematical statement: $X(s,\chi)$ is a diagnostic quantity which directly constrains the  location of their non-trivial zeros. 

\begin{figure}[t]
\centering
\includegraphics[width=0.45\textwidth]{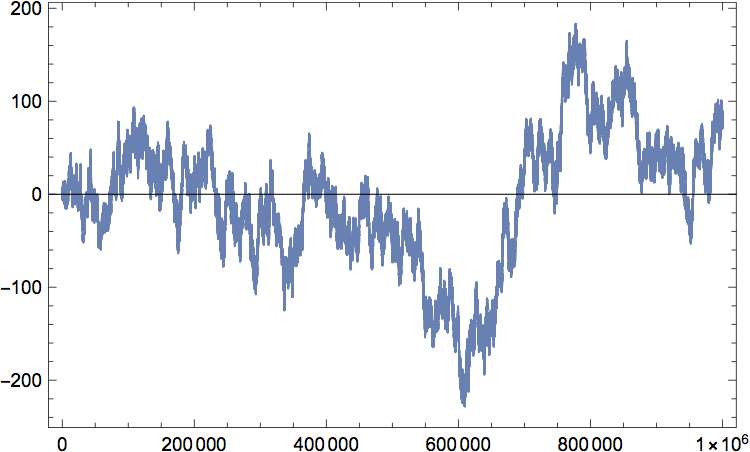}
\caption{Typical plot of $C_N$ versus $N$ (here relative to the character $\chi_2$ mod 7, see Table \ref{table:nonlin}),  
with the same qualitative behavior for any other non-principal characters.}
\label{tipical}
\end{figure}

In fact, since the domain of convergence of  series  like $X(s,\chi)$ are known to be half-planes,  with  $s = \sigma + i t$, we are allowed to consider the case $t=0$. Taking  now the real part\footnote{Analogous arguments apply to the imaginary part of $X(s, \chi)$.}
of $X(s, \chi)$ in \eqref{PDir}, we have (for $t=0$) 
\begin{equation}
\label{SDef}
S(\sigma, \chi) \,=\, \sum_{n=1}^\infty \dfrac{\cos \theta_{p_n}}{p_n^{\, \sigma}} \,\,\,. 
\end{equation}
Defining
\begin{equation}
\label{CxDef}
C(x) \,=\, \sum_{p \le x} \cos \theta_p\,\,\,,  
\end{equation}
we have  
\beq
 C(p_n) - C(p_{n-1}) \,=\, \cos \theta_{p_n}\,\,\,,
\eeq
and then  
\begin{equation}
S(\sigma, \chi) = \sum_{n=1}^{\infty} C(p_n) \left( \dfrac{1}{p_n^{\,\sigma}} - 
\dfrac{1}{p_{n+1}^{\,\sigma}}\right) = \sigma \sum_{n=1}^{\infty} 
C(p_n) \int_{p_n}^{p_{n+1}} \dfrac{1}{u^{\sigma+1}} du\,\,\,.
\end{equation}
Since $C(x) = C(p_n)$ is a constant for $x \in (p_n, p_{n+1})$, we finally arrive to 
\begin{equation}
\label{importantintegral}
S(\sigma,\chi) = \sigma \int_{2}^{\infty} \dfrac{C(x)}{x^{\sigma+1}} \,dx\,\,\,. 
\end{equation}
Hence, the convergence of the integral is determined by the behavior of the function $C(x)$ at $x \rightarrow \infty$: if $C(x) = O(x^{\alpha})$ for $x \rightarrow \infty$, then the integral converges for $ \sigma > \alpha$ and diverges precisely at $\sigma =\alpha$.  Notice that $C(x)$ can change sign (and it does) infinitely many times
(see for instance Figure \ref{tipical}), so in writing $C(x) = O(x^{\alpha})$ we simply mean the order of the function regardless its sign. The conclusion above can be summarised in the following theorem:

\vspace{3mm}
\begin{theorem}\label{seriesinzero}
(Fran\c ca-LeClair \cite{EPFchi})
{\it Consider the sum on the primes 
\beq
\label{Gseries}
C_N \,=\, \sum_{n=1}^{N}    \cos  \theta_{p_n}   
\,\,\,,
\eeq
and assume that such a quantity for large $N$ scales as 
\beq 
\label{scalinglawCN}
C_N  \simeq N^\alpha \,\,\,, 
\eeq
up to logarithms (see the discussion below). Then the GRH for $L$-functions of non-principal characters is true if 
\beq
\alpha = \half  + \epsilon\,\,\,.
\eeq
 for any arbitrarily small $\epsilon > 0$.}
\end{theorem}
Indeed, if $\alpha = \half$ the convergence of the infinite product of the $L$-function can be safely extended down  arbitrarily close to  the critical line  $\Re (s) = \half$ without encountering any zeros.  It is important that extra logarithmic  factors in the scaling law (\ref{scalinglawCN}) of $C_N$, such as $\log N$ or $\log \log N$ etc., do not spoil the convergence argument. For instance, assuming $C(x) = O\left(\sqrt{x} \log^a x \right)$  yields
\begin{equation}
|S(\sigma, t,  \chi)| \le  K \int_{1}^\infty 
\dfrac{\log^a x}{x^{\sigma+1/2}} dx = 
K \dfrac{\Gamma(a+1)}{(\sigma - 1/2)^{a+1}} \,\,\,.
\end{equation}
Henceforth by simply writing $C_N \simeq N^{1/2}$,  without always displaying the $\epsilon$,  implicitly means 
that this scaling behavior can be relaxed with logarithmic factors. 

Since we are interested in the behavior of the series (\ref{Gseries}) only for $N\rightarrow \infty$, notice that we are always free to drop an arbitrarily large number of its first terms. In other words, for the purpose of studying the behavior of $C_N$ in the limit $N \rightarrow \infty$,  one has the following equivalence between series 
\beq
\hspace{25mm}
\sum_{n=1}^N c_n \sim \sum_{n=\tilde N}^N c_n
\,\,\,\,\,\,\,\,\,\,\,\,\,\,
,
\,\,\,\,\,\,\,\,\,\,\,\,\,\,
{\rm for} \, 
N \rightarrow \infty\,\,\,\,\,, 
\label{importantequivalence}
\eeq
where the initial value $\tilde N$ of the second series can be arbitrarily large. 
 
Sums which obey a power law behavior such as $N^{\alpha}$ are common objects in the subject of random walks, and it is useful to recall  that $\alpha = \half $ corresponds to purely {\em diffusive} brownian motion, those with values of $\alpha $ in the interval $ 0 < \alpha < \half$ correspond instead to {\em sub-diffusive} motion, while those in the interval $\half  < \alpha \leq 1$ to {\em super-diffusive} motion, as for instance happens in Levy flights, see \cite{Mazo,Rudnick,Levy,diffusion}. Hence the validity and the universality of the GRH for {\em all} Dirichlet functions of non-principal characters can be linked to a purely diffusive brownian motion behavior\footnote{We can a-priori exclude that the series $C_N$ presents a sub-diffusive behavior since it is known that there are an infinite number of zeros along the critical line $\sigma = 1/2$,  see \cite{Conrey2}.} of the quantity $C_N$. As shown below, such a diffusive behavior of the series $C_N$ can be established on the basis of the Dirichlet theorem on the equi-distribution of reduced residue classes modulo $q$ \cite{Diric,SelbergD} and the LOS conjectures   on the distribution of pairs of residues on consecutive primes \cite{OliverSoundararajan}. For a typical behavior of the series $C_N$, see Figure \ref{tipical}. 

A final remark: in light of the result just presented, it is now easy to see what the problem is with the principal characters $\chi_{1}$: in this case all the angles $\theta_{p_n}$ in the expression of $C_N$ are zero and therefore 
\beq 
C_N   \, =\, \sum_{n=1}^{N} 1 \,=\, N\,\,\,, ~~~~{\rm for} ~ \chi = \chi_1 =  {\rm the ~principal ~character,}
\label{noshiftfunction}
\eeq
i.e. this series grows  {\em linearly} with $N$.  The Mellin transform (\ref{importantintegral}) now diverges at 
$\sigma =1$ but this divergence just signals the pole at $s=1$ of the corresponding $L$-functions and unfortunately provides no information on their zeros.  For this reason,  a proper treatment of this case will not be further considered in this paper, namely here we focus only on the  non-principal  cases.  We should mention that one approach to the principal case involves truncating the Euler product in a well-prescribed manner \cite{Gonek,ALZeta}.

\section{Statistical properties of the angles $\theta_{p_n}$}\label{stprth}

In this Section we are going to discuss some known and conjectured properties of the angles $\theta_{p_n}$ that enter into the series  $C_N$: 
\beq 
C_N \,=\, \sum_{n=1}^{N}    \cos  \theta_{p_n}   \equiv \sum_{n=1 }^{N}  c_n
\,\,\,,
\label{defCNN}
\eeq
The values of these angles belong to a finite and discrete set of angles 
\beq
\theta_{p_n} \in \Phi \,=\, \{ \phi_1,\phi_2,\ldots,\phi_r\}
\,\,\,\,\,\,\,\,
,
\,\,\,\,\,\,\,\,
r \leq \varphi(q) \,\,\,,
\label{setangles}
\eeq
relative to distinct roots of unity. The integer $r$ is commonly referred to as the order of the character and it is  an integer that divides  $\varphi (q)$. Therefore, moving along the sequence of the primes, we are essentially dealing with the rolling of a dice of $r$ faces: the outputs of the various angles, although deterministic, vary in a complicated and irregular way -- too irregular to follow their variation in detail -- but one expects that certain averaged features of this sequence may present a well defined behavior. This is indeed the case and therefore it is perfectly legitimate to enquire about the distribution of the various angles, the frequency of each of them along the sequence of primes, and if the \textquotedblleft dice\textquotedblright \, is biased or not, namely if the various outputs are correlated and how much they are correlated. In other words, in nailing down the behavior of the series $C_N$, we are going to explore the interplay between randomness and determinism, according to a well established and successful tradition in Number Theory (see, for instance, \cite{Dyson,Montgomery,Odlyzko,Rudnick-Sarnack,Conrey,EPFchi,
Kac,Billingsley,GrosswaldSchnitzer,Schroeder,Tao,Chernoff,Torquato,Cramer}). As shown below, this approach has far reaching consequences. 
\def\hatA{A} 

Let us define the sequence of $N$ angles relative to the first $N$ primes: 
\beq
\hatA_{N}  \,=\, 
\{\theta_{p_n}\, ; ~  n = 1, 2, \ldots, N \} \,\,\,. 
\label{sequenceSNpart}
\eeq
As anticipated above, there is important statistical information on the sequence of the angles $\hatA_N $. 
The first important result is the Dirichlet theorem.  As the infinite series of odd numbers $1, 3, 5, \ldots (2 m+1)$ contains infinitely many primes, an interesting question to settle is whether this property is also shared by other arithmetic progressions such as 
\beq
S_m \,=\, q \, m + a \,\,\,\,\,,\,\,\,\, m=0, 1, 2, \ldots 
\,\,\,\,\,\,\,\,\,\,\,\, 
q, a  \in \mathbb{N} 
\label{sequence}
\eeq
In such progressions, the number $q$ is known as the {\em modulus} while the number $a$ as the {\em residue}. It is easy to see that a necessary condition to find a prime among the values of $S_m$ is that the two natural numbers $q$ and $a$ have no common divisors, namely they are {\em coprime}, a condition expressed as $(q,a) = 1$. In 1837 Dirichlet proved that this condition is also sufficient, that is if $(q,a) = 1$ then  the sequence  $S_m$ contains infinitely many primes. His ingenious proof involves some identities satisfied by the Dirichlet $L$-functions. One basic result is the prime number theorem for arithmetic progressions.  Namely,  define
\beq
\label{pntq}
\pi_a (x; q )  = \# \{ p < x\, : p = a ~ {\rm mod} ~ q \}\,\,\,,
\eeq
and 
\beq
\label{nprimes}
\pi (x )  = \# \{ p < x\,\}\,\,\,.
\eeq
Then, for $x \rightarrow \infty$, Dirichlet proved that 
\beq
\label{pntq2}
\pi_a (x; q ) \sim \pi(x)/\varphi (q)\,\,\,.
\eeq
Since  
\beq
\pi(x) \sim {\rm Li}(x)
\,\,\,\,\,\,\,\,\,\,
,
\,\,\,\,\,\,\,\,\,\,
x \rightarrow \infty 
\label{asymptoticpi}
\eeq
where ${\rm Li}(x) \sim x/\log (x)$ is the log integral function, eq.\,(\ref{pntq2}) can be also written as 
\beq 
\pi_a(x; q) \sim {\rm Li}(x)/\varphi(q) \,\,\,.
\label{asympiaq}
\eeq
Since the angles $\theta_{p}$ are functions of the residue $a$ of the prime $p$ mod $q$, one can thus reframe Dirichlet's theorem into the statement that the angles $\theta_{p_n}$ are {\em equally distributed} among their possible $r$ values:   

\vspace{3mm}
\begin{theorem}\label{Dirichlettheorem} (Dirichlet)  
{\it Let $\chi(p_n) = e^{i \theta_{p_n}} \neq 0$ be a non-principal Dirichlet character modulo $q$ and $\pi(x)$ the number of primes  less than $x$. These distinct roots of unity form a finite and discrete set, i.e. $\theta_{p_n} \in \Phi = \{ \phi_1,\phi_2,\ldots,\phi_r\}$ with $r \leq \varphi(q)$ and we have
\beq 
f_i = \lim_{x\rightarrow \infty} 
\frac{\#\{ p\leq x \, : \theta_p = \phi_i\}}{\pi(x)}
= \frac{1}{r} 
\label{equipp}
\eeq
for all $i = 1, 2, \ldots, r$ where $p$ denotes a prime while $f_i$ denotes the frequency of the event $\theta_p = \phi_i$ occurring.
}
\end{theorem}

\vspace{3mm}
\noindent
{\bf Examples.} Inspection of the characters mod $7$  in Table \ref{tablech} indicates that $\chi_4$ has $r=2$, $\chi_{3}$ and $\chi_5 $ have  $r=3$, while $\chi_{2}$ and $\chi_6$ have $r=\varphi (7) = 6$.

\vspace{3mm}
\noindent
{\bf Correlations}. The Dirichlet theorem tells us that, in the limit $N \rightarrow \infty$, the angles $\theta_{p_n}$ are equally distributed in the sequence $\hatA_N $,  but it does not specify, in detail, {\em how} they are distributed. Important insights on the correlations of these variables comes from the knowledge of how many times the pairs $(\phi_a,\phi_b)$ appear as values of two consecutive angles $\theta_{p_n}$ and $\theta_{p_{n+1}}$, or of angles separated by $k$ steps in the sequence $A_{N} (\s)$, i.e. $\theta_{p_n} $ and $\theta_{p_{n+k}}$. 

This problem has been recently addressed by Lemke Oliver and Soundararajan on the basis of the Hardy-Littlewood prime k-tuples conjecture. The basic concept that motivated these works is the so-called Chebyshev bias (see \cite{GranvilleMartin}): for instance, studying prime numbers less than $x$, for {\em finite} values of $x$, there are more primes of the form $4m+3$ than $4m+1$, whereas one naively expects equal numbers. In the paper \cite{OliverSoundararajan}, instead of the angles $\theta_{p_n}$, Lemke Oliver and Soundararajan were directly concerned with the patterns of residues mod $q$ among the sequences of consecutive primes less than an integer $x$ (on this subject see also \cite{Shiu,Ash} and references therein). For our purposes, this is equivalent to the correlations among the angles $\theta_{p_n}$ since these quantities are just functions on the residues.  We will refer to the residues as ``$a$" (or ``$b$")  in accordance to  \eqref{sequence}:  
\beq
\label{resid}
p_n = a ~ {\rm mod} ~ q, ~~~a \in \{ 0,1,2, \ldots , \varphi(q)  \} 
\eeq
If $q$ is not a prime, then not all values of $a$ in the above set are realized. If $q$ is instead a prime, then 
for $p_n > q$,  there are only $\varphi (q)$ possible values of the residue $a$ and only when $p=q$ is the residue equal to $0$. So, except for the prime $q$,  the number of possible residues is $\varphi (q)$.  Henceforth we limit our attention to $q$ equal to a prime,  where $\varphi (q) = q-1$,  but will continue to refer to it as $\varphi (q)$. 
Quantities of interest are the counting functions 
\beq
\pi_{ab} (x; q, k) \,= \#\{ p_n < x  : 
p_n \equiv a\,\,({\rm mod}\,\,q)\,\,,\,\, p_{n+k} \equiv b \,\,({\rm mod}\,\,q) \}
\eeq
For instance,  for $q=3$ and $k=1$,  $\pi_{ab}$ counts the number of consecutive primes whose residues have the patterns $(a,b) = (1,1), (1,2), (2,1), (2,2)$. Based on the pseudo-randomness of the primes, for $x \rightarrow \infty$ one would expect that  the primes counted by $\pi_{ab}(x; q, k)$ go as 
\beq 
\pi_{ab}(x;q,k) \sim \pi(x)/(\varphi(q))^2\,\,\,,
\label{largeasymptotic}
\eeq
independent of the separation $k$ of the two residues but, as noticed by  Lemke Oliver and Soundararajan, for {\em finite values of $x$} there are large secondary terms in the expected asymptotic behavior which create biases toward and against certain patterns of residues. In the following, in particular, we focus our attention on the $\varphi (q) \times  \varphi (q) $ matrices\footnote{LOS define $f_{ab}$ as above but with $\pi (x)$ replaced by the log integral ${\rm Li} (x)$.  The latter is simply the leading approximation to $\pi (x)$ based on the prime number theorem, thus our definition is actually more meaningful. In the large $x$ limit,  the results are the same whether one uses $\pi (x)$ or ${\rm Li }(x)$.}
\beq
f_{ab}(x,q,k) \,=\, \frac{\#\{p_n \leq x : p_n \equiv a \,\,({\rm mod}\,\,q)\,\,,\,\, p_{n+k} \equiv b \,\,({\rm mod}\,\,q) 
\}}{\pi(x)} \,\,\,, 
\label{definitionfreqOliver}
\eeq
which, for nearby $x$, give the {\em local densities} of pairs of primes, in which $p_n \equiv  a \,\,({\rm mod}\,\,q) $ will be followed, after $k$ steps, by a prime $p_{n+k} \equiv b \,\,({\rm mod}\,\,q)$. Here we quote the large $x$ behavior of these quantities \cite{OliverSoundararajan}: 

\vspace{1mm}

\begin{conjecture}(Lemke Oliver-Soundararajan). 
{\it For large values of $x$ we have\footnote{It is possible to express $f_{ab}(x,q,1)$ and $f_{ba}(x,q,1)$ individually (and they are not equal) but their expression is rather complicated, see \cite{OliverSoundararajan}. Moreover, the expressions (\ref{firstconjab}) and (\ref{firstconjaa}) given here are those of LOS  but specialised to the modulus $q$ being a prime.} 
\beq 
f_{ab}(x,q,1) + f_{ba}(x,q,1) \,=\,\frac{2}{(\varphi(q))^2} \left[1 + \frac{\log\log x}{2 \log x} - \log\frac{q}{2\pi} 
\frac{1}{2 \log x} + O\left(\frac{1}{(\log x)^{7/4}}\right)\right]\,\,\,,
\label{firstconjab}
\eeq
whereas 
\beq
f_{aa}(x,q,1) \,=\, \frac{1}{(\varphi(q))^2} \left[1 - \frac{(\varphi(q) -1)}{2} \, \frac{\log\log x}{\log x} + 
(\varphi(q) -1) \,\log\frac{q}{2\pi} \,\frac{1}{2 \log x} +
O\left(\frac{1}{(\log x)^{7/4}}\right)\right]\,\,\,.
\label{firstconjaa}
\eeq
} 
\end{conjecture}

\vspace{1mm}

\begin{conjecture}(Lemke Oliver-Soundararajan). 
{\it  For  $k \geq 2$,  then for large values of $x$ we have 
\beq 
f_{ab}(x,q,k)  \,=\,\frac{1}{(\varphi(q))^2} \left(1 + \frac{1}{2 (k-1)} \,\frac{1}{\log x} + O\left(\frac{1}{(\log x)^{7/4}}\right)\right)\,\,\,,
\label{secconjab}
\eeq
whereas 
\beq
f_{aa}(x,q,k) \,=\, \frac{1}{(\varphi(q))^2} \left(1 - \frac{(\varphi(q) -1)}{2(k-1)} \, \frac{1}{\log x} + O\left(\frac{1}{(\log x)^{7/4}}\right)\right)\,\,\,.
\label{secconjaa}
\eeq
} 
\end{conjecture}
 
\vspace{1.5mm}
\noindent
{\bf Remark.} The opposite signs in \eqref{firstconjaa} versus \eqref{secconjab} are responsible for the bias conjectured by LOS. Moreover, notice that these correlations have a permutation symmetry ${\mathcal S}_{\varphi(q)}$ since the only thing that matters is whether the residues are equal or different.

\vspace{3mm}
\noindent

We can use these formulas to study the correlations among the angles in the subsequences $\hatA_N$.
In particular,  the counting functions $ \pi_{ab}(x;q,k)$ of the pairs of primes in $\hatA_N $ relative to various residues is given by
\beq
\label{generaldistrrr}
 \pi_{ab}(x;q,k)  =  \
\pi (x) \, f_{ab}(x,q,k) \, \sim \frac{x}{\log x}  f_{ab} (x, q, k)
\eeq
The interesting aspects of these functions $\pi_{ab}(x;q,k)$ are the following: 
\begin{enumerate}

\item For $x\to  \infty$, all pairs of residues in $\hatA_N  $ are {\em equally probable} (both for consecutive primes and primes separated by $k$ steps) and their probability is given by $1/(\varphi(q))^2$. This means that, in the limit $x \rightarrow \infty$, the angles $\theta_{p_n}$ in any subsequence $\hatA_N $ are completely {\em uncorrelated}. 
\item However, at any {\em finite} value of  $x$, the next neighbor variables in the subsequence $\hatA_N $ tend to be {\em anti-correlated}: the occurrence of pairs of equal residues $(a,a)$ for next neighbor primes are always less probable than the occurrence of pairs of different residues $(a,b)$. In other words, there is a bias in the distribution of the residues  between consecutive primes: notice, however, that once again this is a finite-size effect that vanishes as $\sim \log\log x/\log x$.  
\item At any {\em finite} $x$, this anti-correlation phenomenon also persists for primes which are separated by $k$ steps and the matrices $f(x,q,k)$ are {\em not} equal to $(f(x,q,1))^k$, i.e. these probabilities do not satisfy the Markovian property. This correlation decreases as $1/k$ with the separation $k$ of the two primes but it is also a finite size effect since the coefficient in front of this $1/k$ correlation vanishes as $1/\log x$ when $x \rightarrow \infty$. The Markovian property of these matrices is of course restored in the $x \rightarrow \infty$ limit.

\end{enumerate}

\def\chat{\hat{c}}
\def\Chat{\hat{C}}
\def\Ghat{\hat{G}}

\def\hatA{A}

\section{Stochastic Time series}\label{Timeseries}

In this section, we bring results of the last section to bear on our mail goal, which concerns the growth of the series $C_N$ as a function of $N$ for large $N$.  Ref.\,\cite{longpaper} also provides extensive numerical 
evidence for the results below.

\subsection{The Single Brownian Trajectory Problem}

Let's recall that our aim is to estimate how the series $C_N$, defined in eq.\,(\ref{defCNN}), grows with $N$. 
We have already remarked that this series presents stochastic features and probabilistic aspects that justify regarding  it as a random time series. Still, in order to fully exploit this point of view, we have to face the problem of defining an ensemble $\CE$ of the possible outputs of $C_N$, together with their relative probabilities. A-priori this issue seems to pose a severe obstacle to this approach since, for any given character of the $L$ functions, there is {\em one and only one} series $C_N$ to deal with. This, however, is a common problem in many time series, in particular for all those that refer to situations for which is impossible to 
\textquotedblleft turn back time\textquotedblright : indeed, in these cases it is impossible to have access to all possible outputs of the relative variable and therefore equally impossible to define the relative probabilities. In the literature, this is known as the {\em Single Brownian Trajectory Problem} (see, for istance \cite{brow1,brow2,brow3,brow4} and references therein). 

A way to get around this problem is to consider an arbitrarily long time series and take 
\textquotedblleft stroboscopic\textquotedblright \,snapshots of it. 
We will do this in the following specific manner.    Define the  {\it ordered} intervals of length $N$ starting at $\s$  
\beq
I_N(\s) =\{\s, \s+1, \s+2, \ldots, \s+N -1 \} \,\,\,, 
\,\,\,
\label{intervals}
\eeq 
and  the associated angles $\hatA_{N}(\s) $ 
\beq
\hatA_{N} (\s)  \,=\, 
\{\theta_{p_n}\, ; ~ n \in I_N(\s)\}\,\,\,. 
\label{sequenceSNpart}
\eeq
We then define {\it block variables} $C_N (\s)$ based on the above intervals:
\beq
C_N(\s) \,=\,\sum_{k\in I_N(\s)}  c_{k} = \sum_{k=\s}^{\s + N -1}  \cos \theta_{p_n}  \,\,\,. 
\label{groupvariables}
\eeq
Notice that, with this new notation, our previous series $C_N$, given in eq.\,(\ref{defCNN}), is simply expressed as $C_N(1)$.  For a greater flexibility in making some of the arguments below, let us also define
\beq 
A_{N_1,N_2} \,  \equiv   A_{N_2 - N_1 +1} (N_1)  = \, \{\theta_{p_n}\, ; ~ n=N_1, N_1+1,\ldots , N_2 \}\,\,\,,
\label{sequenceSN}
\eeq
relative to primes between $p_{N_1}$ and $p_{N_2}$. 
When 
\beq
N_1 < \s < \s + N -1 < N_2\,\,\,\,,
\label{appartenenza}
\eeq
$\hatA_{N}(\s)$ is of course a subset  of $A_{N_1,N_2}$.  Let's clarify the role of $N_1$ and $N_2$. Imagine we fix a very large value of $N_1$ and then vary $N_2$: in this way we can consider arbitrarily long sequences $A_{N_1,N_2}$, out of which many and well separated block variables $C_N(\s)$ of the {\em same length} $N$ can be extracted and used as members of the ensemble to which belongs the original sequence $C_N(1)$! This is equivalent to the stroboscopic snapshots behind the solution of the the Single Brownian Trajectory Problem (see the forthcoming subsection): the validity of this self-averaging procedure relies on two aspects of the corresponding time series, its ergodicity and stationarity. 

In the case of our sequences $A_{N_1,N_2}$, their ergodicity is guaranteed by {\em all} possible outputs of the angles $\theta_{p_n}$ along the sequence of the primes. Their stationarity is an issue more subtle which nevertheless can be settled on the basis of the following considerations: according to the formulas of LOS , there are correlations which explicitly depend on the point $x$ along the sequence of the primes and therefore, for arbitrary values of the extrema $N_1$ and $N_2$, they break -- strictly speaking -- the stationarity of the sequences $A_{N_1,N_2}$. 

There are however two facts which help in solving this issue: the first is that, as we already commented, these are finite size effects which vanish when $x\rightarrow \infty$; the second is the equivalence (\ref{importantequivalence}) which, even at finite $x$, makes us free to focus our attention on sequences whose extrema $N_1$ and $N_2$ are such that the correlations are both weak and sufficiently uniform along the entire length of these intervals. Intervals $(N_1,N_2)$ which satisfy this property will be called {\em inertial intervals} and sequences based on these intervals can be made stationary as much as one wishes. For instance, choosing $N_1 = 10^{200}$ and $N_2 = 10^{250}$, the correction to a uniform background $1/(\varphi(q))^2$ distribution is only of the order $0.20 \%$ and $0.17\%$ respectively at the beginning and at the end of the sequence  $A_{N_1,N_2}$, therefore with a breaking of the stationarity that can be quantified of the order of $0.03 \%$. These values come from the correction $1/\log x$ present in the LOS with respect to  the constant values of the correlations, computed for $x =N_1$ and $x=N_2$. If one is unhappy with this $0.03\%$ breaking of the stationarity, one can choose instead, say,  $N_1 = 10^{1000}$ and $N_2 = 10^{1200}$ so that these corrections at the beginning and at the end of the interval drop to the more negligible values $0.043 \%$ and $0.036\%$, with a breaking of the stationarity of the series of only  $0.007\%$. 

Given that the number of primes is infinite, the point is that we can always choose higher and higher values of $N_1$ and $N_2$ and make the corresponding sequence $A_{N_1,N_2}$ stationary with any arbitrary degree of confidence. By the same token, namely enlarging at our wish the size of the sequences $A_{N_1,N_2}$, we can always set up a proper ensemble for $C_N(1)$ for {\em any} $N$, no matter how large. Notice that as $N_{1}$ and $N_2 \to \infty$, also $\s \to \infty$. Moreover, we are going to assume the inequalities  
\beq
\label{ineq}
1 \ll N \ll \s,  
\eeq
so that  $p_\s \approx p_{\s+N}$.    

\subsection{Statistical Ensemble $\CE$ for  the series $C_N$}\label{subesemble}
\begin{figure}[t]
\centering
\includegraphics[width=0.40\textwidth]{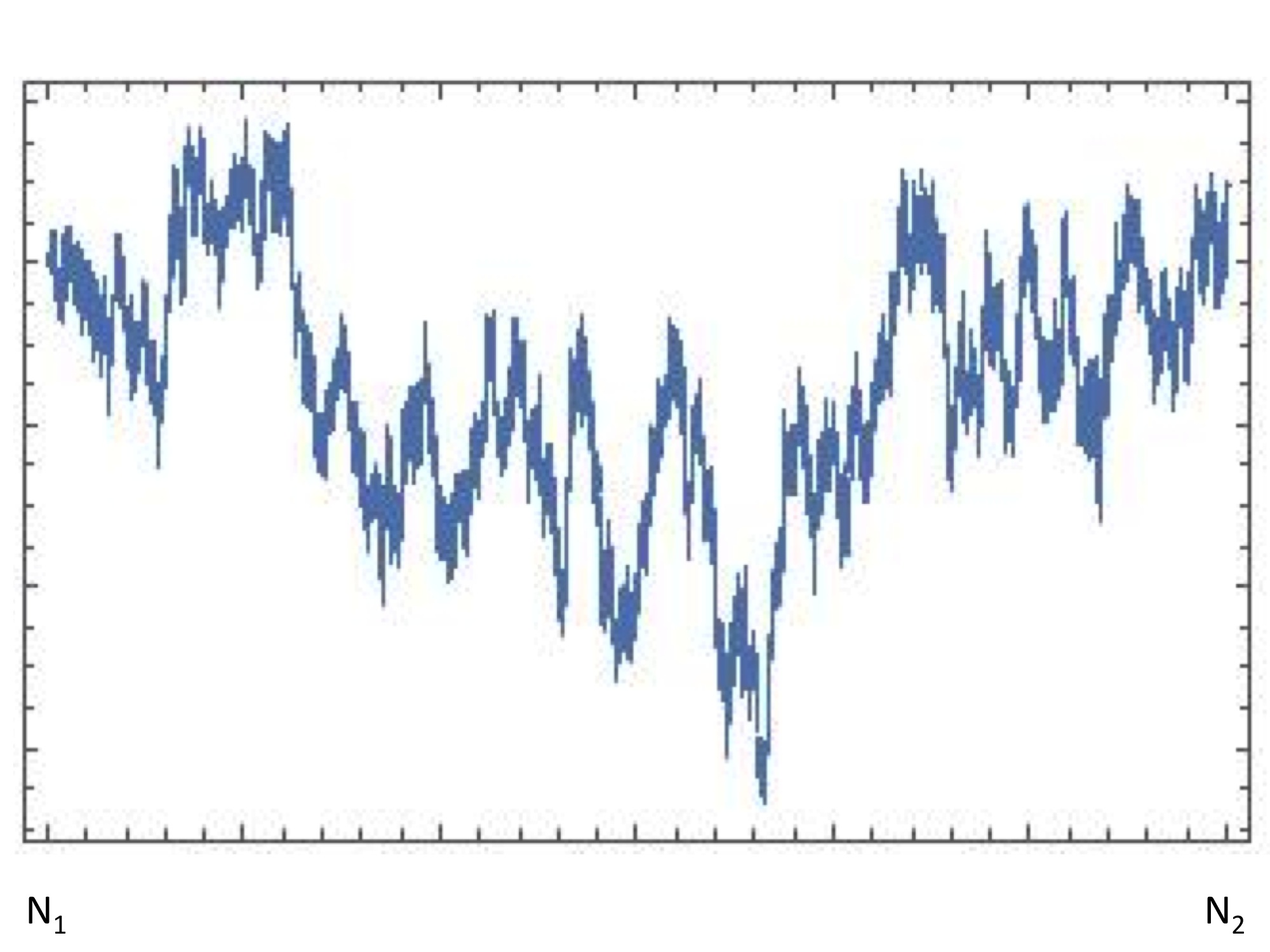}
\,\,\,
\includegraphics[width=0.40\textwidth]{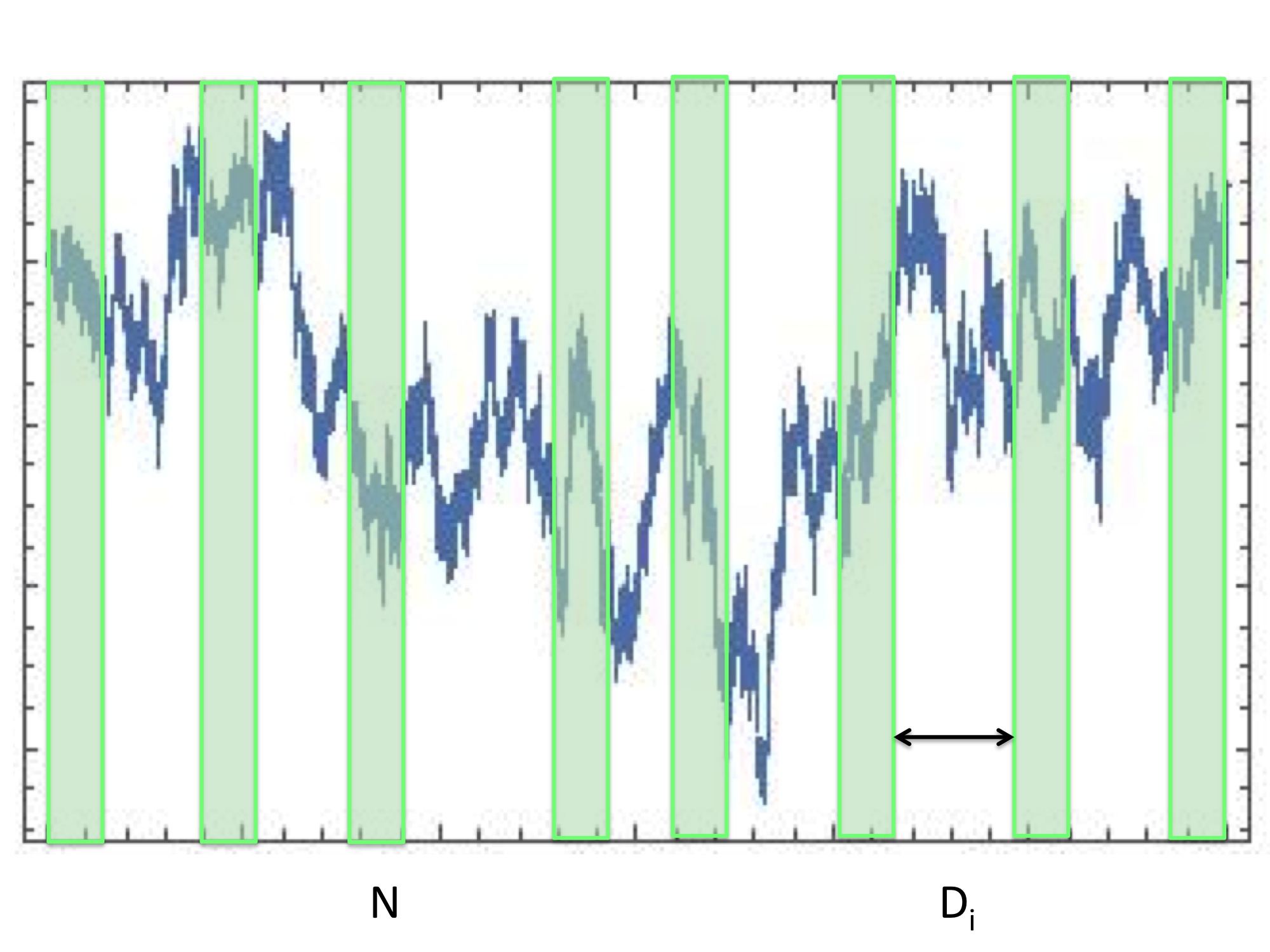}
\caption{Left hand side: series $\sum_{n=N_1}^{N} $ vs $N$, in the inertial interval $(N_1,N_2)$. 
Right hand side: sampling of the time series done in terms of block variables $C_N (\s)$ of length $N \gg 1$ relative to the green intervals separated by distances $D_i$. Under the hypothesis of stationarity of the sequence $A_{N_1,N_2}$, the values of these blocks define a probabilistic ensemble $\CE$ for the quantity $C_N$ relative to the sum of the first $N$ values $c_n$. } 
\label{partitiontimeseriesGN}
\end{figure}

The block variables $C_N(\s)$ are the equivalent of the  
\textquotedblleft stroboscopic\textquotedblright \,images of length $N$ of a single Brownian trajectory (see Figure \ref{partitiontimeseriesGN}) and they allow us to control the irregular behavior of the original series $C_N(1)$ by proliferating it into a collection of sums of the same length $N$. It is this collection of sums that forms the {\em set of events}, i.e. the ensemble $\CE$ relative to the sums of $N$ consecutive terms $c_n$. The procedure to set up such an ensemble is as follows: 

\begin{enumerate}
\item Consider two very large integers $N_1$ and $N_2$  (which eventually we will send to infinity), with $N_1\gg 1$, $N_2\gg1$  but also $L \equiv (N_2 - N_1) \gg 1$ such that, for a given character $\chi$ of modulus $q$,  
the sequence $A_{N_1,N_2}$ is inertial. 
\item For any fixed integer $N$, with $1 \ll N \ll L$, consider the sets 
\beq 
{\cal S}_M \,=\, \bigcup_{i=1,\ldots M} I_N(i) 
\,\,\,\,\,\,\,\,
,
\,\,\,\,\,\,\,\,
N_1 \leq i < N_2 
\,\,\,\,\,\,,\,\,\,\,\, I_N(i) \cap I_N(j) = 0 
\,\,\,,\,\,\,  i \neq j 
\eeq
made of $M$ {\em non-overlapping} and also well separated intervals of length $N$ whose origin is between the two large numbers $N_1$ and $N_2$ (see Figure \ref{partitiontimeseriesGN}). These conditions ensure that the block variables $C_N(i)$ computed on such disjoint intervals are very weakly correlated and therefore we can assume that we are dealing statistically with $M$ separated copies of the original series $C_N(1)$. 
\item At any given $N_1$ and $N_2$, the cardinality ${\rm card} ({\cal S}_M) = M$ of these sets cannot be larger of course than $L/N$. There is however a large freedom in generating  them: 
\begin{enumerate}[label={\bf \alph*}]
\item We can take, for instance, $M$ intervals $I_N(\s)$ separated by a fixed distance $D$, with the condition that $M(N + D) = L$; 
\item Alternatively,  we can take,  $M$ intervals $I_N(\s)$ separated by random distances $D_i$ such that $M N + \sum_{i=1}^M D_i = L$. 
\end{enumerate}
\item The ensemble $\CE$ is then defined as  the set of the $M$ block variables $C_N(\s)$ relative to the intervals $I_N(\s) \in {\cal S}_M$:
\beq
\label{CEdef}
\CE = \{ C_N (\s) \},  ~~~~{\rm with~~}  I_N (\s) \in {\cal S}_M
\eeq
\end{enumerate}

\noindent 
In summary, choosing two very large and well separated  integers $N_1$ and $N_2$, we can generate a large number of sets of intervals ${\cal S}_M$ and use the corresponding block variables of length $N$ to sample the typical values taken by a series consisting of a sum of $N$ consecutive terms $c_n$. In view of the ergodicity and stationarity of the sequence $A_{N_1,N_2}$ for $N_1\rightarrow \infty$ and $N_2\rightarrow \infty$, this is tantamount to determining  the statistical properties of the original series $C_N$. 

\subsection{Mean and variance of the block variables $C_N(\s)$}
Let's now use the ensemble $\CE$ to compute the most important quantities of the series $C_N(\s)$, i.e. the mean and the variance of these sums. The value of the mean is a simple consequence of the Dirichlet theorem.

\begin{figure}[t]
\centering\includegraphics[width=.6\textwidth]{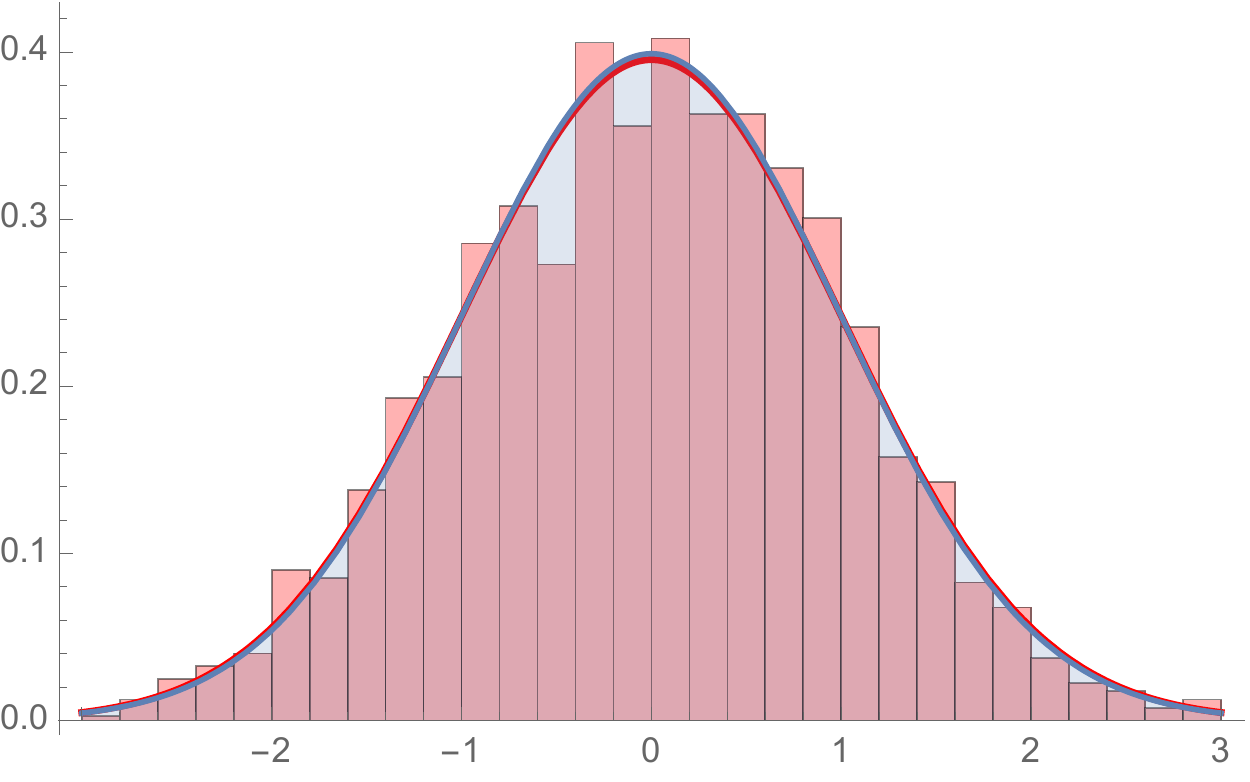}
\caption{Numerical evidence for the normal distribution proposed in \eqref{finallynormal}.     What is shown is a histogram of the LHS of  \eqref{finallynormal}   which are properly normalized block variables $C_N  (\s)$ for the character $\chi_2$ mod $7$ in Table \ref{tablech}.    The ensemble $\CE$ corresponds to $N=10,000$,  $D=10$,  with 
$M=2000$ states. The red curve is the  fit to the data,  which is the normal distribution $\CN(0.00026, 1.009)$.   The  nearly indistinguishable  blue curve is the prediction $\CN(0,1)$.}
\label{FigureCLT1}
\end{figure}

\vspace{3mm}
\noindent
{\bf Mean of $C_N$}. In the limit $N\rightarrow \infty$, the series $C_N$ has zero mean 
\beq 
\mu\, \equiv \, \lim_{N\rightarrow \infty} \frac{1}{N} 
\sum_{n=1}^N \cos\theta_{p_n} \,=\, 0 \,\,\,.
\label{meanC_N}
\eeq
The proof is immediate. Consider the case when the cardinality $r$ of the set $\Phi$ of the angles coincides with $\varphi(q)$, i.e. $r = \varphi(q)$. We can use then eq.\,(\ref{pairwiseangles}) to group pairwise the terms of the sum and 
since  
\beq
\cos(\alpha_{\varphi/2 +k}) = - \cos\alpha_{k}
\,\,\,\,\,\,
,\,\,\,\, 
k=1,\ldots ,\varphi(q)/2 \,\,\,,  
\label{pairwisecosine}
\eeq
we have 
\beq 
\mu \,=\,\, \lim_{N\rightarrow \infty} \frac{1}{N} 
\sum_{n=1}^N \cos\theta_{p_n} \,=\, 
 \sum_{k=1}^{\varphi(q)/2} \cos\theta_k \,
 \left(f_k - f_{\varphi/2+k}\right)
 \,=\, 0\,\,\,
\label{meanC_N11}
\eeq
since, in the $N\rightarrow \infty$ limit, from  the Dirichlet theorem all frequencies $f_{n}$ are equal. Analogous results can be easily obtained also when $r < \varphi(q)$. In the double limit $N_1 \rightarrow \infty$ and $N_2 \rightarrow \infty$ (so that also $N \rightarrow \infty$), from the stationarity properties of the sequence $A_{N_1,N_2}$ the same is true for the ensemble average of the large $N$ block variables $C_N(\s)$ 
\beq 
\Ex \Bigl[C_N  (\s )\Bigr] \,=\, 0 \,\,\,.
\eeq
In light of this result, the ensemble $\CE$ consists of block variables $C_N(\s)$ which are equally distributed among positive and negative values. The typical histogram of these block variables $C_N(\s)$ assumes the bell-shape shown in Figure \ref{FigureCLT1}.     Our next step is to compute how the variance of the block variables scales as function of $N$.  

\vspace{3mm}
\noindent
{\bf Variance of the block variables $C_N(\s)$}. The block variables $C_N(\s)$ are defined in eq.\,(\ref{groupvariables}).
Let us  first define the variance $b^2$ of the cosine on the set of the $r$ angles 
\beq
\label{variancecosine}
b^2 \,\equiv\,\frac{1}{r} \sum_{k=1}^r \cos^2\phi_k \,=\, 
\left\{
\begin{array}{lll} 
1 & \,\,\,\,,\,\,\,\, & {\rm \,\,}   {\rm if ~ \chi ~ is ~ real} \\
1/2 &  \,\,\,\,,\,\,\,\, & {\rm  \,\,} {\rm if ~ \chi ~ is ~ complex } 
\end{array}
\right.
\eeq 
If  $\chi$ is real,  then the only values of the character are $\chi = \pm 1$. Of course, if the terms $c_n$ of these sums were uncorrelated, it would have been easy to compute the probability distribution of these block variables in terms of the characteristic function $\hat P(k)$ of the variable $c\equiv \cos\theta$. The formula would have been simply given by  
\beq 
P(x) \,=\,\frac{1}{2\pi} \int_{-\infty}^{\infty} dk (\hat P(k))^N \, e^{-i k x} \,\,\,, 
\label{PX}
\eeq
and this would have led immediately to the gaussian behavior relative to the central limit theorem, since 
\beq 
\hat P(k) \,\simeq 1 - b^2\,\frac{k^2}{2} + \cdots \,\,\,,
\eeq
with $b^2$ given in eq.\,(\ref{variancecosine}), and for large $N$ 
\beq
P(x) \,=\,\frac{1}{2\pi}\,\int_{-\infty}^{\infty} dk \,e^{N \log \hat P(k)} \, e^{-i k x} \,
\simeq \frac{1}{2\pi}\,\int_{-\infty}^{\infty} dk \,e^{-N b^2 \,k^2/2}\ \, e^{-i k x} \,
\simeq e^{-x^2/(2 N b^2)}\,\,\,. 
\eeq
In other words, if the $c_n$'s were uncorrelated, the block variables $C_N(\s)$ would have been straightaway gaussian distributed with a variance equal to $N$ times the variance $b^2$ of the $c_n$'s. 

However, for any sequence $A_{N_1, N_2}$, the variables $c_n$ {\em are} correlated (although weakly) and the actual computation of the variance $\sigma_N^2$ of the block variables $C_N(\s)$ must be done using the formulas of LOS. Such a computation goes as follows\footnote{Here we present the argument relative to the case $r=\varphi(q)$ but the final expression of the variance, eq.\,(\ref{mostimportantformula}), holds for all cases. Moreover, in the following we will consider block variables nearby the position ${\bf \s}$, with $N_1 \leq \s < N_2$. 
}.  
Consider the block variable $C_N(\s)$ belonging to the ensemble $\CE$ defined above and take the ensemble average of its square 
\begin{eqnarray}
\sigma_N^2  (\s) & \,=\, & \Ex \Bigl[ (C_N(\s))^2 \Bigr] \,=\,\, \sum_{l=0}^{N-1}  \sum_{m=0}^{N-1} 
\Ex \Bigl[ c_{\s + l } c_{\s + m } \Bigr]\,=\,\sum_{m=0}^{N-1} \Ex \Bigl[ c^2_{\s+m}\Bigr] + 
2 \, \sum_{m=1}^{N-1} [N - m] \,\Ex \Bigl[ c_{\s} c_{\s+ m} \Bigr] \,\,\,\nonumber \\
&\,=\,&\, N \, \Ex \Bigl[ c^2_\s \Bigr] + 2 \, \sum_{m=0}^{N-1} (N - m) \Ex \Bigl[ c_{\s } c_{\s + m} \Bigr] \,\,\,,
\end{eqnarray}
where we used the stationarity of the ensemble to group the contributions of the various pairs separated by $k$ steps (there are $(N-m)$ of them). Isolating further the term $m=1$ in the second quantity of the expression above, we have that the variance can be expressed as 
\beq
\Ex \Bigl[ (C_N(\s))^2 \Bigr] \,=\, {\cal D}_0 + {\cal D}_1 + {\cal D}_2 
\eeq
where 
\begin{eqnarray}
{\cal D}_0 \,& = & \,  N \, \Ex \Bigl[ c^2_\s\Bigr] \,\,\,,\nonumber \\
{\cal D}_1 \,& = & \,  2  (N - 1) \,\Ex \Bigl[ c_{\s} c_{\s + 1} \Bigr] \,\,\,,\\
{\cal D}_2 \,& = & \, 2 \sum_{m=2}^{N-1}  [N - m] \,\Ex \Bigl[ c_{\s} c_{\s + m} \Bigr] \,\,\,.\nonumber
\end{eqnarray}
The variables $c^2$ are statistically equi-distributed on the $r$ angles and their variance $b^2$ on these angles was given in eq.\,(\ref{variancecosine}), so ${\cal D}_0$ is expressed as  
\beq
{\cal D}_0 \,=\,b^2 \,  N \,\,\,. 
\label{d00}
\eeq
In order to compute ${\cal D}_1$, we need the formulas (\ref{firstconjab}) and (\ref{firstconjaa}) relative to the residues of two next neighbor primes. In light of eqs.\,(\ref{firstconjab}) and (\ref{pairwisecosine}), notice that in  the average of the product of the two cosines on the ensemble, keeping initially the $\cos\theta_{p_i}$ fixed, there are only {\em two terms} which contribute to the average: the first when $\theta_{p_{i+1}}  \,=\,\theta_{p_i}$ (with weight $f_{aa}(p_\s,q,1)$), the second when $\theta_{p_{i+1}} \,=\,\theta_{p_i} + \pi$ (with weight $f_{ab}(p_\s,q,1)$), while all other terms cancel out pairwise. Summing now on the $\varphi(q)$ values taken by $\theta_{p_i}$, we have 
\beq
{\cal D}_1 \,=\, - b^2\, (N-1)  \,\frac{\log\log p_\s}{\log p_\s} +b^2\, \log\frac{q}{2\pi} \frac{1}{\log p_\s}\,\,\,. 
\label{d11}
\eeq

The calculation is essentially similar for the other term ${\cal D}_2$, the only difference between the dependence of the separation $m$ of the two cosines
\beq
{\cal D}_2 \,=\, - b^2 \,\frac{1}{\log p_\s}\,\sum_{m=2}^{N-1} (N - m) \frac{1}{m-1}
\,\,\,.
\eeq
Putting together the three terms, we arrive to the following theorem: 
\vspace{3mm}
\begin{theorem}\label{GM1}   
{\it Assuming the validity of the LOS  conjectures, in the inertial intervals 
the variance $\sigma_N^2$ of the the block variables of length $N$ is given by  
\beq
\sigma_N^2(\s)/b^2\,=\, \Ex \Bigl[ (C_N(\s))^2 \Bigr]/b^2 \,=\,  N \,\lambda(N,\s) \,  + \,\rho(N,\s) \,\,\,,
\label{mostimportantformula}
\eeq
where 
\begin{eqnarray}
\lambda(N,\s) &\,=\,& \left[1 +\frac{1}{\log p_\s} \left(1 - \sum_{m=1}^{N-2} \frac{1}{m}\right) - \frac{\log\log p_\s}{\log p_\s}\right] \,\,\,,
\label{correctionfactorlambda}\\
\rho(N,\s) &\,=\,& \, \frac{1}{\log p_\s} \left[\log \left(\frac{q \log p_\s}{2 \pi e^2}\right) + \sum_{m=1}^{N-2} \frac{1}{m}
\right] 
\label{interceptrho}\,\,\,.
\end{eqnarray}
}
\end{theorem}

\vspace{3mm}
Theorem \ref{GM1} is our main result: in all the inertial intervals, the variance of block variables $C_N(\s)$ of length $N$ scales {\em linearly} with N, up to a correction factor $\lambda(N,\s)$ which is {\em independent} of the modulus $q$,  but depends on the prime $p_{\s}$ of the inertial interval around which we consider the block variables. Notice that, for $\s\rightarrow \infty$, we recover a purely gaussian expression for  the variance 
\beq
\lim_{\s \rightarrow \infty} \sigma^2_N (\s)  \,=\, b^2 \, N \,\,\,.  
\label{limitxinfinity}
\eeq
Note that we have used the inequality \eqref{ineq} which implies $p_{N+\s} \approx p_\s$. 
Keeping instead $\s$ finite and considering the large $N$ asymptotic of this formula, the factor $\lambda(N,\s)$ introduces a logarithmic correction since 
\beq
\sum_{m=1}^{N-1} \frac{1}{m} \simeq \log N + \gamma_{E} \,\,\,,
\eeq
where $\gamma_E$ is the Euler-Mascheroni constant. Notice that, as far as $\s$ is finite, for the anti-correlation of the residues of consecutive primes, we have $\lambda(N,\s) < 1$ and therefore the variance of the block variables $C_N(\s)$ at a finite $\s$ is always {\em smaller} than the variance of $N$ uncorrelated variables.

Theorem \ref{GM1},  along with  \eqref{limitxinfinity},  implies that in the limit $\s \to \infty$ 
the properly normalized block variables are gaussian distributed: 
\beq
\frac{C_N(\s)}{\sigma_N(\s)} ~ \dist ~  \CN(0,1)\,\,\,.
\label{finallynormal}
\eeq
where finite $\s$  corrections to  $\sigma_N(\s)$  are given in eq.\,(\ref{mostimportantformula}).  For finite $\s$, the distribution is not purely gaussian, and this non-gaussianity captures the existence of correlations between the primes for a given character.  
Numerical evidence for this normal distribution is shown in Figure  \ref{FigureCLT1}.

\vspace{3mm}

In light of this result, taking larger and larger inertial intervals, and therefore correspondingly larger and larger values of $N$,  the block variables $C_N(\s)$ of length $N$ always scale as 
\beq
C_N(\s) \,\simeq {\cal O}(N^{1/2+\epsilon}) \,\,\,. 
\label{fffff}
\eeq
for arbitrarily small $\epsilon >0$.  
Notice that, in probabilistic language, in the limit $N \rightarrow \infty$ this behavior occurs  with probability equal to 1. Indeed, since $\sigma_N(\s) \leq \sqrt{N}$, using the normal law distribution (\ref{finallynormal}), in the limit $N\rightarrow\infty$ we have 
 \barray
 \nonumber 
 \prob \[  |C_N(\s)  |    <   d  \, \sqrt{N} \]  & > &
 \prob \[  |C_N(\s)   |    <   d  \,  \sigma_N(\s)  \]   \,=\,  
 \inv{\sqrt{2 \pi}} \int_{-d}^d  dx\,  e^{-x^2/2}  
 \\ 
 \label{PrBN} 
 &=& 1 -  \frac{e^{- d^2 /2}} {\sqrt{ 2 \pi } } \(  \frac{2}{d}  +  O\(  \inv{d^2 }  \)  \).
 \earray
  Chose $d  = \kappa N^\epsilon$ for any $\kappa >0$.   
 Then for any $\epsilon > 0 $,   
 \beq
 \label{ProbOne}
 \lim_{N \to \infty}  \prob \[  C_N(\s)  = O(N^{1/2+ \epsilon}) \]  = 1    \,\,\,.
 \eeq

However, we can show that the actual result is even stronger than this probabilistic argument, namely   
we can explicitly exclude the existence of  ``rare events"  that would potentially  spoil the purely diffusive behavior of the series $C_N$. Let's use a reductio ad absurdum argument, namely let's assume that, in view of some rare events\footnote{For our series this would be for example a very long series of the same residues for successive primes.}, the series $C_N$ for $N \rightarrow \infty$ rather than going as in eq.\,(\ref{fffff}) would instead behave  as $C_N \simeq N^{\alpha}$ (up to logarithmic corrections) with $\alpha \neq \half$. It this were true, such a behavior of the series $C_N$ should hold for any neighbourhood of infinity, namely for all $N$ which satisfies 
$N > N^*$, for any arbitrarily large $N^*$. But, in turn, this would imply that the variance of $C_N$ should always go as $N^{\alpha}$ for {\em all} the infinitely many ensembles $\CE$ of the inertial sequences $A_{N_1,N_2}$ with $N_1 > N^*$. The infinite occurrence of such behavior would contradict firstly the notion itself of \textquotedblleft rare events\textquotedblright \, and, secondly, it would be in clear contrast with the explicit expression (\ref{mostimportantformula}) of the variance computed on all the infinitely many ensembles $\CE$ of the inertial intervals. In other words, we cannot exclude that {\em some} $C_N$ for some specific starting point $\tilde \s$ of the series (\ref{importantequivalence}) and even for long values of $N$ may grow as $N^\alpha$ with $\alpha \geq  \half$, but, if this would be the case, using the equivalence of the various series related to $C_N$,  we can always change at our will $\tilde \s$ as well as we can also take larger and larger values of $N$. Choosing the new $\tilde \s$ to be inside any of the inertial intervals, a behaviour of the block variable as $N^\alpha$ would then disappear in favor  of the only stable behavior of the series $C_N$ under any possible translation of the inertial intervals and any possible ensemble $\CE$ set up in these intervals, namely the scaling law of the random walk given by $N^{1/2}$ (again up to logarithms).   

\section{Conclusions}\label{conclusions}

In this paper we have addressed the Generalised Riemann Hypothesis for the Dirichlet $L$-functions of non-principal characters. A diagnostic quantity for the location of all non-trivial zeros of these functions is given by the large $N$ behavior of the series $C_N$ defined in eq.\,(\ref{defCNN}): a purely diffusive random walk behavior as 
$N^{1/2 +\epsilon}$ of this series, for arbitrarily small $\epsilon >0$, signifies that all zeros of these functions are along the critical line $\Re(s) = \half$. We have shown that there is a natural explanation of such a diffusive behavior of the series $C_N$, in light of the following circumstances: 
\begin{enumerate}
\item The terms of this series, given by $c_n = \cos\theta_{p_n}$, involve the angles of the characters $\chi(p_n)$ computed on the primes $p_n$. Therefore,  what  matters are the residues of the primes with respect a chosen modulus $q$ rather than the primes themselves, and this turns the problem into studying the outputs of a  \textquotedblleft dice\textquotedblright \, of $\varphi(q)$ faces. 
\item There is  certain important information on the distribution of these angles $\theta_{p_n}$. Their equi-distribution is guaranteed by the Dirichlet theorem on arithmetic progressions while their correlations are ruled by the conjectures of LOS. By virtue of the Dirichlet theorem, $C_N$ has, on  average, as many positive as negative terms of equal values. The absence in $C_N$ of one (or more terms) much larger than all others makes a-priori difficult, if not impossible, for the series to have a superdiffusive behavior $N^{\alpha}$ with $\alpha > \half$, as for instance happens in Levy flight where the motion is characterised by clusters of shorter jumps inter-sparsed by long jumps. In the case of $C_N$, instead, any significant growth of this quantity must necessarily be the result of many small increments all of the same sign, an event which is probabilistically extremely rare, of the order $1/(\varphi(q))^m$, where $m$ is the number of these equal outputs.  Since we are only interested in the behavior of $C_N$ at $N \rightarrow \infty$, all these fluctuations average out as in the usual random walk problem and lead to a purely diffusive behavior. There are some correlations among consecutive angles $\theta$'s which however do {\em not} spoil the diffusive behavior of the series $C_N$ coming from the equidistribution of its positive and negative values. 
\item As in the {\em single Brownian trajectory problem}, in order to extract the properties of the series $C_N$ we have taken the point of view to consider the inertial sequences $A_{N_1,N_2}$ of the angles $\theta_{p_n}$ (see eq.\,\ref{sequenceSN})) as time series.  In view of the ergodicity and stationarity of the $A_{N_1,N_2}$, we have set up an ensemble $\CE$ for the series $C_N$ in terms of the block variables $C_N(\s)$ defined on well separated subintervals of length $N$ with origin at $\s$ along these inertial sequences. 
\item Based on Dirichlet theorem and assuming the LOS conjecture, we have shown that block variables of length $N$ satisfy a normal distribution, with a variance which goes linearly in $N$ (up to logarithmic corrections) as shown in eq.\,(\ref{mostimportantformula}). 
\end{enumerate}

Hence, based on Theorem \ref{GM1}, which assumes the LOS conjectures, we have established  a purely diffusive random walk behavior of the series $C_N$. The most important consequence concerns the diagnosis of the non-trivial zeros of the Dirichlet $L$-functions of non-principal characters, related to the divergence of the integral (\ref{importantintegral}). The conclusion is that this integral diverges at $\Re(s)  = \half$ and this implies 
that all zeros  of the Dirichlet $L$-functions of non-principal characters all of them are along the same critical line.    A natural question is how strongly this conclusion depends on the LOS conjectures? We would answer that 
the most important property the pair correlation is their asymptotic uncorrelated behavior given in eq.\,(\ref{largeasymptotic}), while the details of the LOS formula are essential for controlling finite $\s$ effects.

\vspace{5mm}
 
\section*{Acknowledgments}

GM would like to thank Don Zagier, Karma Dajani, Andrea Gambassi and Satya Majumdar for interesting discussions and Robert Lemke Oliver for useful email correspondence.  AL would like to thank SISSA in Trieste, Italy, where this work was begun while GM would like to thank the Simons Center in Stony Brook and the 
International Institute of Physics in Natal for the warm hospitality and support during the initial and final parts of this work respectively. 


\vspace{3mm}

\appendix

\section{Dirichlet Characters}\label{Appcharacters}
Given a modulus $q$, the prime residue classes modulo $q$ form an abelian group, denoted as   
\beq
(\mathbb{Z}/q\mathbb{Z})^* := \{m\, {\rm mod} \,q \,:\, (m, q) = 1\} \,\,\,. 
\label{groupabbb}
\eeq
The dimension of this group is given by the Euler totient arithmetic function $\varphi(q)$ which counts 
how many positive integers less than $q$ are coprime to $q$.  A Dirichlet character $\chi$ of modulus $q$ is an arithmetic function from the finite abelian group $(\mathbb{Z}/q\mathbb{Z})^*$ onto $\mathbb{C}$ satisfying the following properties: 
\begin{enumerate}
\item 
$\chi(m+q) \,=\, \chi(m) $. 
\item 
$\chi(1) =1 $ and $\chi(0) = 0$. 
\item 
$\chi( m \,n ) \,=\, \chi(m) \, \chi(n)$. 
\item 
$\chi(m) = 0$ if $(m,q) > 1$ and $\chi(m) \neq 0$ if $(m,q) =1$. 
\item 
If $(m,q) =1$ then $(\chi(m))^{\varphi(q)} =1$, namely $\chi(m)$ have to be $\varphi(q)$-roots of unity. 
\item 
If $\chi$ is a Dirichlet character so is the complex conjugate $\overline\chi$. 
\end{enumerate}
From property $5$, it follows that for a given modulus $q$ there are $\varphi(q)$ distinct Dirichlet characters that can be labeled as $\chi_{j}$ where $j = 1, 2, . . . , \varphi(q)$ denotes an arbitrary ordering (we will not always display this index $j$ in $\chi_j$). Moreover the characters satisfy the following orthogonality conditions 
\begin{eqnarray}
\sum_{r=1}^{\varphi(q)} \chi_r(k) \overline{\chi}_r(l) &\,=\, & 
\left\{ 
\begin{array}{cll}
\varphi(q) & & {\rm if} \,\,k\equiv l \,\,\, ({\rm mod} \, q) \\
0 & & {\rm if} \,\, k\not\equiv l \,\,\, ({\rm mod} \, q) 
\end{array}
\right. \\
\sum_{m=1}^{q} \chi_r(m) \overline{\chi}_s(m) & \,=\,& 
\varphi(q) \, \delta_{r,s} \,\,\,.
\end{eqnarray}
For a generic $q$, the {\em principal} character, usually denoted $\chi_1$, is defined as 
\beq
\chi_1(m) \,=\, 
\left\{ 
\begin{array}{cl}
 1 & \, {\rm if} \,\,(m, q) = 1 \\
 0 & \, \, {\rm otherwise} 
 \end{array}
 \right.
 \eeq
When $q = 1$, we have only the {\em trivial}  principal character $\chi(m) = 1$ for every $m$, and in this case the corresponding $L$-function reduces to the Riemann $\zeta$-function given by 
\beq
\zeta(s) \,=\,\sum_{m=1}^\infty \frac{1}{m^s} 
\,\,\,, 
\,\,\,\,\,\,\,\,
\mathbb{R}(s) > 1 \,\,\,. 
\label{zetafunctions}
\eeq 

There is an important difference between principal versus non-principal characters. The principal characters, being only $1$ or $0$, satisfy 
\beq
\sum_{m=1}^{q-1} \chi_1(m) \,=\,\varphi(q) \,\neq \,0 \,\,\,,
\eeq
whereas  the non-principal characters  satisfy
\beq
\sum_{m=1}^{q-1} \chi(m) \,=\,0 \,\,\,.
\label{sumtozero}
\eeq

\vspace{3mm}
\noindent
{\bf Parametrization of the angles}. Posing 
\beq
\label{thetan2}
\chi (m) \,=\,  e^{i \theta_m }, ~~~~~\forall ~\chi(m) \neq 0 \,\,\,, 
\eeq   
eq.\,(\ref{sumtozero}) shows that the angles $\theta_m$ of the non-principal characters defined in eq.\,(\ref{thetan}) are equally spaced over the unit circle, and being associated to the $\varphi(q)$ roots of unity, their 
possible values can be labelled as 
\beq
\label{notationangles}
\alpha_k \,=\, \frac{\pi (2 k - \varphi(q))}{\varphi(q)} \,\,\,\,\,\, , 
\,\,\,\,\,\, k = 1,\ldots, 
\varphi(q) \,\,\,.
\eeq 
In this parameterization the angles $\alpha_k$ are negative for $k=1,\ldots\,\varphi(q)/2$ while they are positive for $k=\varphi(q)/2+1,\ldots,\varphi(q)$, related pairwise as 
\beq
\alpha_{\varphi/2 +k} \,=\,\alpha_k + \pi \,\,\,\,\,\,
,\,\,\,\, 
k=1,\ldots ,\varphi(q)/2 \,\,\,.
\label{pairwiseangles}
\eeq
Notice that the actual distinct roots of unity entering the expression of the characters may be a smaller set of the $\varphi(q)$-roots of unity, $\theta_m \in \Phi = \{ \phi_1,\phi_2,\ldots,\phi_r\}$ with $r \leq \varphi(q)$ and  
$\phi_i$ equal to one of the angles of the set (\ref{notationangles}). The integer $r$ is referred to as the {\em order} of the particular character. 

As an explicit example of the various characters associated to a modulus $q$, consider $q=7$ where they are expressed in terms of the $6$-th roots of unity, as shown in Table \ref{tablech}, with $\omega = e^{i\pi/3}$.  Here $\varphi(7)=6$.  
Notice that $\chi_1$ and $\chi_4$ are real (and the corresponding angles belong to a smaller set of the $6$-roots of unity) while the terms of the pairs $(\chi_2,\chi_6)$ and $(\chi_3,\chi_5)$ are complex conjugates of each other. The characters $(\chi_3,\chi_5)$ are composed  of  an $r=3$ subset of the angles (\ref{notationangles}) (i.e. the angles $0, \pm 2\pi/3$), whereas  the characters $(\chi_2,\chi_6)$ employ the full set of angles (\ref{notationangles}). 

\begin{table}[t]\label{tablech}
\centering
\begin{minipage}{.4\textwidth}
\centering
\begin{tabular}{| c | c | c | c | c | c | c | c |}
\hline\hline
n & 1 & 2 & 3 & 4 & 5 & 6 & 7 \\ [0.5ex] 
\hline
$\chi_1(n)$ 
& \,\, 1 \,\,& 1 & 1 & 1 & 1 & \,\,1 \,\,& \,\,\, 0 \,\,\\
\hline
$\chi_2(n)$ & 1 & $\omega^2$ & $\omega$ & $-\omega$ & $-\omega^2$ & \,\,-1 \,\,& \,\,\, 0 
\,\,\\
\hline
$\chi_3(n)$ & 1 & $- \omega$ & $\omega^2$ & $\omega^2$ & $-\omega$ & \,\, 1 \,\,& \,\,\, 0 \,\,\\
\hline
$\chi_4(n)$ & 1 & 1 & -1 & 1 & -1 & \,\,-1 \,\,& \,\,\, 0 \,\,\\
\hline
$\chi_5(n)$ & 1 & $\omega^2$ & $-\omega$ & $-\omega$ & $\omega^2$ &\,\, 1 \,\, & \,\,\, 0 \,\,\\
\hline
$\chi_6(n)$ & 1 & $-\omega$ & $-\omega^2$ & $\omega^2$ & $\omega$ & \,\,-1 \,\, & \,\,\, 0 \,\, \\
\hline
\end{tabular}
\end{minipage}
\caption{Characters mod $q=7$, where 
$\omega = e^{i \pi/3}$. }
\label{table:nonlin}
\end{table}

\section{Functional Equation for the Dirichlet $L$ functions}
Here we present  the functional equation satisfied by the Dirichlet $L$-functions, similar to the one of the Riemann $\zeta$-function (for the proof of this equation, see for instance \cite{Apostol}). To express such a functional equation, let's define  $a$ for  a character as  
\beq\label{order}
a \equiv \begin{cases}1 \qquad
\mbox{if $\chi(-1)= -1$ ~~\,\,\,\,(odd)} \\
0 \qquad \mbox{if $\chi(-1) = \,\,\,\,\,1$  ~~~~(even).}
\end{cases}
\eeq
Let's also introduce the Gauss sum 
\beq
G(\chi) \,=\,\sum_{m=1}^q \chi(m) \, e^{2 \pi i m/q}
\,\,\,. 
\label{Gausssum}
\eeq
$G(\chi)$ satisfies $|G(\chi)|^2 = q$ if and only if the character $\chi$ is primitive. With these definitions, the functional equation for the primitive $L$-functions can be written as 
\beq
\label{FELambda}
L(1-s, \chi) \,=\, i^{-a} \,\frac{q^{s-1} \, \Gamma(s)}{(2\pi)^s} 
\, G(\chi) \, \left\{\begin{array}{c}\cos(\pi s/2)\\ \sin(\pi s/2) 
\end{array}\right\} \, 
L(s, \overline\chi)\,\,\,. 
\eeq
where the choice of cosine or sine depends upon the sign of $\chi(-1) = \pm 1$.  

\section{Residue at the pole $s=1$} 
As emphasized,  there is an important distinction between the $L$-functions based on non-principal versus principal characters, 
and this is their residue at $s=1$. To show this result, one first expresses  any $L$-function in terms of a {\em finite} linear combination of the Hurwitz zeta function $\zeta(s,a)$ as 
\beq
L(s,\chi) \,=\, \frac{1}{q^s} \, \sum_{r=1}^q \chi(r) \, \zeta\left(s, \frac{r}{q}\right) 
\,\,\,.\nonumber
\label{LHurwitz}
\eeq 
where 
\beq
\zeta(s,a) \,=\, \sum_{m=0}^\infty \frac{1}{(m+a)^s} \,\,\,,
\label{Hurwitz}
\eeq
The Hurwitz $\zeta$-function has a simple pole at $s=1$ with residue 1 and therefore the residue at this pole of the $L$-function is 
\beq
{\rm Res} \, L(s,\chi) \,=\, \frac{1}{q} \sum_{r=0}^q \chi(r)\,=\, 
\left\{
\begin{array}{cll}
\frac{\varphi(q)}{q} & & {\rm if }\, \, \chi = \chi_1 \\
0 & & {\rm if }\,\, \chi \neq \chi_1 \,\,\,.
\end{array}
\right.
\eeq
Thus, 
\begin{itemize}
\item
The $L$ functions for  non-principal characters are {\em entire} functions,  i.e. analytic everywhere in the complex plane with no poles. 
\item The $L$-functions  $L(s,\chi_1)$ for  principal characters, on the contrary, are analytic everywhere except for a {\em simple pole} at $s=1$ with residue $\varphi(q)/q$. 
\end{itemize}

\end{document}